\documentclass[10pt]{article}

\usepackage{amsmath,latexsym,amsthm,amsfonts,amssymb}

\newtheorem{theorem}{Theorem}

\newtheorem{corollary}{Corollary}
\newtheorem{observation}{Observation}
\newtheorem{proposition}{Proposition}
\newtheorem{definition}{Definition}

\newcommand{\IN}{{\mathbb N}}

\newcommand{\w}{\omega}

\newcommand{\Op}{\mathcal O}
\newcommand{\U}{\mathcal U}
\newcommand{\A}{\mathcal A}
\newcommand{\B}{\mathcal B}
\newcommand{\PP}{\mathcal P}
\newcommand{\J}{\mathcal J}
\newcommand{\F}{\mathcal F}
\newcommand{\G}{\mathcal G}
\newcommand{\Fr}{\mathfrak{F}r}

\newcommand{\utas}{\bigcup_{\mathrm{fin}}(\mathcal O,\mathrm{T}^\ast)}
\newcommand{\uta}{\bigcup_{\mathrm{fin}}(\mathcal O,\mathrm{T})}
\newcommand{\uop}{\bigcup_{\mathrm{fin}}(\mathcal O,\mathcal O)}
\newcommand{\uga}{\bigcup_{\mathrm{fin}}(\mathcal O,\Gamma)}
\newcommand{\uom}{\bigcup_{\mathrm{fin}}(\mathcal O,\Omega)}
\newcommand{\utass}{\bigcup_{\mathrm{fin}}(\mathcal O,\mathrm{T}^\star)}
\newcommand{\dep}{\mathrm{Depth}^+([\omega]^{\aleph_0})}

\newcommand{\muu}{\mathfrak u}
\newcommand{\mb}{\mathfrak b}
\newcommand{\md}{\mathfrak d}
\newcommand{\ms}{\mathfrak s}
\newcommand{\mg}{\mathfrak g}
\newcommand{\mt}{\mathfrak t}

\title{A semifilter approach to selection principles II:
$\tau^\ast$-covers}
\author{Lyubomyr Zdomskyy}
\begin{document}
\maketitle

\begin{abstract}
Developing the ideas of \cite{Zd} we show that every Menger topological space
has the property $\utas$ provided $(\mathfrak u<\mathfrak g)$, 
and every space with the property $\utas$
is Hurewicz provided $(\dep\leq\mb)$. Combining this with the results proven in cited literature, we settle all
 questions whether (it is consistent that) the properties $\mathsf P$ and $\mathsf Q$
[do not] coincide, where $\mathsf P$ and $\mathsf Q$ run over
$\uga$, $\uta$, $\utas$, $\uom$, and $\uop$.  
\end{abstract}

\large
\centerline{\textbf{Introduction}}
\normalsize
\medskip

Following \cite{Sc} we say that a topological space
$X$ has the property $\bigcup_{\mathrm{fin}}(\mathcal A,\mathcal B)$, where $\mathcal A$ and
$\mathcal B$ are collections of covers of $X$, if for every sequence $(u_n)_{n\in\w}\in\mathcal A^\w$
there exists a sequence $(v_n)_{n\in\w}$, where each $v_n$ is a finite subset of
$u_n$, such that $\{\cup v_n:n\in\w\}\in\mathcal B$. Throughout this paper ``cover'' means
``open cover'' and 
$\mathcal A$ is equal to the family $\Op$  of all open covers of $X$.
Concerning $\mathcal B$, we shall also consider the collections $\Gamma$, $\mathrm{T}$, $\mathrm{T}^\star$, $\mathrm{T}^\ast$, 
and $\Omega$
of all open $\gamma$-, $\tau$-, $\tau^\star$, $\tau^\ast$-, and $\w$-covers   of $X$.
For technical reasons we shall  use the collection $\Lambda$ of countable large covers. 
The most natural way to define these types of covers uses the Marczewski ``dictionary'' map
introduced in \cite{M-S}.  
\footnotetext{\normalsize \emph{Keywords and phrases.} Selection principle,
 semifilter,  small cardinals. 

\emph{2000 MSC.} 03A, 03E17, 03E35.\footnotesize }
Given  an indexed family $u=\{U_n:n\in\w\}$  
of subsets of a set $X$ and element $x\in X$,  we define the Marczewski map $\mu_u:X \to \PP(\w)$
letting $\mu_u(x)=\{n\in\w:x\in U_n\}$ ($\mu_u(x)$ is nothing else but $I_s(x,u)$ in notations of \cite{Zd}).
Recall, that  $A\subset^\ast B$ means that $|A\setminus B|<\aleph_0$.  A family
$\A\subset\PP(X)$ of  subsets of a set $X$ is a \emph{refinement} of a family $\B\subset\PP(X)$,
 if for every $B\in \B$ there exists $A\in\A$ such that $A\subset B$. 
Depending on the properties of $\mu_u(X)$ a family $u=\{U_n:n\in\w\}$ is defined to be
\begin{itemize}
\item a \emph{large cover} of $X$ \cite{Sc}, if for every $x\in X$ the set $\mu_u(x)$ is infinite;
\item a \emph{$\gamma$-cover} of $X$ \cite{GN}, if for every $x\in X$ the set $\mu_u(x)$
is cofinite in $\w$, i.e.  $\w\setminus \mu_u(x)$ is finite;
\item a \emph{$\tau$-cover} of $X$ \cite{Ts-tow},  if it is a large cover and 
the family $\mu_u(X)$ is linearly
preordered by the almost inclusion relation $\subset^\ast$
in sense that for all $x_1,x_2\in X$ either $\mu_u(x_1)\subset^\ast \mu_u(x_2)$
or $\mu_u(x_2)\subset^\ast \mu_u(x_1)$; 
\item a \emph{$\tau^\ast$-cover} of $X$ \cite{Ts-tow}, if  there exists 
a  linearly preordered by $\subset^\ast$
refinement $\J$ of $\mu_u(X)$ consisting of infinite subsets of $\w$; 
\item an \emph{$\w$-cover} \cite{GN}, if the family $\mu_u(X)$ is centered,
i.e. for every finite subset $K$ of $X$ the intersection $\bigcap_{x\in K}\mu_u(x)$
is infinite. 
\end{itemize}
We also introduce a new type of covers situated between $\tau$- and $\tau^\ast$-covers.
A family $u=\{U_n:n\in\w\}$ is 
\begin{itemize}
\item a \emph{$\tau^\star$-cover} of $X$, if there exists a linearly preordered by $\subset^\ast$
refinement of $\J\subset\mu_u(X)$ of $\mu_u(X)$
consisting of infinite subsets of $\w$. 
\end{itemize} 
Recall, that  $\uga$  and $\uop$
are nothing else but the well-known Hurewicz and Menger covering properties introduced
in \cite{Hu} and \cite{Me} respectively at the beginning of 20-th century.

Since every $\gamma$-cover is a  $\tau$-cover, every $\tau$-cover is a $\tau^\star$-cover,
every $\tau^\star$-cover is  a $\tau^\ast$-cover, and every $\tau^\ast$-cover is an $\w$-cover, 
the above properties are related as follows:
\medskip

\begin{center}
\small
\begin{tabular}{ccccccc}
 $\uta$       & $\Rightarrow$ & $\utass$   & $\Rightarrow$ & $\utas$    & $\Rightarrow$ & $\uom$       \\
$\ ^{(2)}$    &               & $\ ^{(3)}$ &               & $\ ^{(4)}$ &               & $\ ^{(5)}$  \\              
 $\Uparrow$   &               &            &               &            &               & $\Downarrow$ \\                               
$\uga$        &               &            &               &            &               &  $\uop$      \\          
$\ ^{(1)}$    &               &            &               &            &               &  $\ ^{(6)}$          

\end{tabular}
\normalsize
\end{center}
\medskip

By a \emph{tower} we understand a $\subset^\ast$-decreasing transfinite
sequence of infinite subsets of $\w$, i.e. a sequence $(T_\alpha)_{\alpha<\lambda}$
such that $T_\alpha\subset^\ast T_\beta$ for all $\alpha\geq\beta$.
The cardinality $\lambda$ is called the \emph{length} of this tower. 
The subsequent theorem, which is the main result of this paper, describes when
some of the above properties coincide. 

\begin{theorem} \label{main}
\begin{itemize}
\item[$(1)$] Under $(\muu<\mg)$ the selection principles $\utass$ and \\
$\uop$ coincide.
\item[$(2)$] Under Filter Dichotomy the selection principles
$\utass$ and $\uom$ coincide.
\item[$(3)$] Selection principles $\uga$ and $\utas$ coincide iff each
 semifilter generated by a tower is meager.
\end{itemize}
\end{theorem}
The following statement describes some partial cases of 
Theorem~\ref{main}(3).

\begin{corollary} \label{details}
\begin{itemize}  
\item[$(1)$] Selection principles $\uga$ and $\utas$ coincide  if
the inequality $\dep\leq\mb$ holds.
\item[$(2)$] Under $(\mb<\md)$ (resp. $(\mt=\md)$) there exists a set of reals with the property
$\utas$  which fails to satisfy $\uga$ (resp. $\uta$).
\end{itemize}
\end{corollary}

Theorem~\ref{main} gives a partial answer of Problem 5.2 from \cite{BST}. 
Namely, it implies the subsequent
\begin{corollary}
It is consistent that the property $\uta$ is closed under unions of families
of subspaces of the Baire space of size $<\mb$.
\end{corollary}
\begin{proof} 
Follows immediately from Theorem~\ref{main}(3) and the fact that 
the property $\uga$ is preserved by unions
of less than $\mb$ subspaces of the Baire space, see \cite{JMSS}.
\end{proof}

We refer the reader to \cite{Va} for definitions of all small cardinals
and related notions  we use. All notions concerning semifilters may be
found in \cite{BZ} and will be defined in the next section.
The condition $(\muu<\mg)$ is known to be consistent: 
 $\mathfrak u=\mathfrak b=\mathfrak s<\mathfrak g=\mathfrak d$ in Miller's model
 and the inequality $(\muu<\mg)$ implies $\muu=\mb<\mg=\md$, see \cite{BlassHBK} and \cite{Va}.
Moreover, $(\muu<\mg)$ is equivalent to an assertion that all upward-closed neither
meager nor comeager families of infinite subsets of $\w$ are ``similar'',
see \cite{La}, \cite[9.22]{BlassHBK}, \cite[7.6.4, 12.2.4]{BZ},
or Theorem~\ref{lll}. This assertion together with the Talagrand's \cite{Ta} characterization
of meager and comeager upward-closed families is a so-called \emph{trichotomy}
for upward-closed families or \emph{Semifilter Trichotomy} in terms of
\cite{BZ}.  The \emph{Filter Dichotomy} follows from the Semifilter Trychotomy
and is formally stronger than the principle NCF introduced by A.~Blass,
see \cite[\S~9]{BlassHBK} and references there in.

 $\dep$ denotes the smallest cardinality $\kappa$ such that there is no
tower of length $\kappa$. Thus $\mt<\dep$.
A model with $\mb\geq\dep$ was constructed in  \cite{dordal}.
Some other applications of $\dep$ in Selection Principles may be found in \cite{ST}.

Theorem~\ref{main} with   results proven in  \cite{JMSS}, 
\cite{Ts-tow}, \cite{TZ}, and \cite{Zd}, enable us to settle almost all questions whether (it is consistent that)
the properties $\mathsf P$
and $\mathsf Q$ [do not] coincide, where $\mathsf P$ and $\mathsf Q$ run
over $\uop,$ $\uom,$ $\utas,$ $\utass,$ $\uta,$ and $\uga$.
(In fact, we settle all of the questions omitting $\utass$.)
Some   sufficient conditions for $\mathsf P=\mathsf Q$
and  $\mathsf P\neq\mathsf Q$  are  summarized in  Table~1.
Each entry $((i),(j))$, $i\neq j$, contains:
\begin{itemize}
\item  A condition  which implies $(i)=(j)$ (resp. $(i)\neq (j)$) 
provided $i<j$ (resp. $i>j$) or ``?'' if  no such a condition is known;
\item ZFC, if $(i)\neq (j)$  in ZFC and $i>j$;
\item $-$, if $(i)\neq (j)$ in ZFC and $i<j$;
\end{itemize}
and a reference to where this is proven. For example,
``$[x]+[y]$, $[z]$'' means that the sufficiency of the corresponding condition
 was proven in $[z]$, and it can be simply derived by combining
results of $[x]$ and $[y]$. Throughout the table, $\lambda$ stands for $\dep$.

\begin{center}
\tiny
%\scriptsize
%\footnotesize 

\begin{tabular}{|c|c|c|c|c|c|c|} \hline
\multicolumn{7}{|c|}{\small \textbf{Table 1} \tiny} \\ \hline
         &             (1)                          &  (2)                                          &       (3)                                  &       (4)          &          (5)                     &   (6)                            \\  \hline 
(1)  &                                          & $(\lambda\leq\mb)$                               &  $(\lambda\leq\mb)$                           & $(\lambda\leq\mb)$    & --                                & --                                          \\
         &                                          &   Cor.~\ref{details}                          &  Cor.~\ref{details}                        & Cor.~\ref{details} & \cite{BZ-mult}, \cite{CP}, \cite{TZ}       &  \cite{BZ-mult},\cite{CP},\cite{TZ},  \\  \hline                    
(2)   & $\mb<\ms$                                &                                               & $(\lambda\leq\mb)$                            & $(\lambda\leq\mb)$    & --                                & --                                          \\
         & \cite{Ts-tow}+\cite{ST}                  &                                               & Cor.~\ref{details}                         & Cor.~\ref{details} & \cite{TZ}                  & \cite{TZ}                    \\ \hline                    
(3) & $(\mb<\ms)\vee (\muu<\mg)$               & $(\muu<\mg      )$                            &                                            &  $(\lambda\leq\mb)$   & Filter Dich.                 & $(\muu<\mg)$                               \\
         & \cite{Ts-tow}+\cite{ST},\cite{TZ}+Th.~1 & \cite{TZ}+Th.~\ref{main}                     &                                            & Cor.~\ref{details} &  Th.~\ref{main}                  & Th.~\ref{main}                     \\ \hline       
(4)  & $(\mt=\md)\vee (\mb<\md)$                & $(\mt=\md)\vee(\muu<\mg)$                     & ?                                          &                    &  Filter Dich.                & $(\muu<\mg)$                               \\
         & Cor.~\ref{details}                       & Cor.~\ref{details}, \cite{TZ}+Th.~\ref{main} &                                            &                    &  Th.~\ref{main}                  & Th.~\ref{main}                      \\ \hline                     
(5)   & ZFC                                      & ZFC                                           & $\lambda\leq\mb$                              & $(\lambda\leq\mb)$                &                                 & $(\muu<\mg)$                               \\
         & \cite{TZ},\cite{CP},\cite{BZ-mult}       & \cite{TZ}                                     & \cite{TZ}+Cor.~\ref{details}               & \cite{TZ}+Cor.~\ref{details} &                                  &  Th.~\ref{main},\cite{Zd}  \\ \hline 
(6)   & ZFC                                      & ZFC                                           &  $(\lambda\leq\mb)\vee \mathrm{CH}$           &  $(\lambda\leq\mb)\vee\mathrm{CH}$ & CH                     &                                           \\ 
         & \cite{TZ},\cite{CP},\cite{BZ-mult}     & \cite{TZ}                                     &  \cite{TZ}+Cor.~\ref{details}, \cite{JMSS} & \cite{TZ}+Cor.~\ref{details}, \cite{JMSS} & \cite{JMSS} &                                     \\  \hline 
\end{tabular}
\end{center}
\normalsize

\medskip

\large
\centerline{\textbf{Semifilters}}
\normalsize
\medskip

Our main tool is the notion of a semifilter. Following \cite{BZ} a family
$\F$ of nonempty subsets of $\w$ is called a \emph{semifilter}, if for every
$F\in \F$ and $A\ ^\ast\supset F$  the set $A$ belongs  to $\F$.
For example, each family $\A$ of infinite subsets of $\w$ generates
the minimal semifilter $\uparrow\A=\{B\subset\w:\exists A\in\A(A\subset^\ast B)\}$ containing $\A$.
The family $\mathsf{SF}$ of all semifilters contains the smallest element
$\Fr$ consisting of all cofinite subsets of $\w$, and the largest one,
$[\w]^{\aleph_0}$, i.e. the family of all infinite subsets of $\w$.
Throughout this paper
by a \emph{filter} we understand a semifilter which is closed under finite
intersections of its elements. 

Since every semifilter $\mathcal F$ on  $\w$ is a subset of the
powerset $\mathcal P(\w)$, which can be identified with the Cantor space $\{0,1\}^\w$,
we can speak about topological properties of semifilters. Recall, that a subset
of a topological space is \emph{meager}, if it is a union of countably many
nowhere dense subsets. The complements of meager subsets are called \emph{comeager}.
We shall often use the subsequent characterization of meagerness of semifilters
due to Talagrand, see \cite{Ta} and \cite[5.3.1]{BZ}.
\begin{theorem} \label{Tal}
A semifilter $\F$ on  $\w$ is  meager
if and only if there exists an increasing number sequence $(k_n)_{n\in\w}$ 
such that every $F\in \F$ meets all but finitely many half-intervals $[k_n,k_{n+1})$.
\end{theorem}  

A crucial role in the proof of Theorem~\ref{main} belongs to the following 
fundamental result
of C.~Laflamme \cite{La}. Following \cite{BZ}, semifilter $\mathcal F$ on $\w$ is said to be
\emph{bi-Baire}, if it is neither meager nor comeager. Note, that 
there is no comeager filter on $\w$, see \cite[5.3.2]{BZ}.
\begin{theorem} \label{lll} The following conditions are equivalent:
\begin{itemize}
\item[$(1)$] $(\mathfrak u<\mathfrak g)$;
\item[$(2)$] for any bi-Baire semifilters $\F$ and $\U$ 
  there exists an increasing
number sequence $(k_n)_{n\in\omega}$  such that 
the sets $\{\{n\in\omega:F\cap [k_n,k_{n+1})\neq\emptyset\}:F\in\mathcal F\}$
and $\{\{n\in\omega:U\cap [k_n,k_{n+1})\neq\emptyset\}:U\in\mathcal U\}$ 
coincide.
\end{itemize}
\end{theorem}

Thus the inequality $(\muu<\mg)$ implies the \emph{Filter Dichotomy} \cite[9.16]{BlassHBK}, which is the abbreviation
of the assertion of Theorem~\ref{lll}(2) for bi-Baire filters:
\medskip

 {\it For arbitrary bi-Baire filters  $\mathcal F$ and $\U$ 
  there exists an increasing number
 sequence $(k_n)_{n\in\omega}$   such that 
the sets $\{\{n\in\omega:F\cap [k_n,k_{n+1})\neq\emptyset\}:F\in\mathcal F\}$
and $\{\{n\in\omega:U\cap [k_n,k_{n+1})\neq\emptyset\}:U\in\mathcal U\}$ 
coincide.}
\medskip

 The main idea of the semifilter
approach to selection principles is to assign to a topological space $X$ 
the family $\{\uparrow\mu_u(X):u\in\Lambda(X)\}$. As it was shown in \cite{Zd},
the  property $\uop$ of a space $X$ may be characterized in terms of topological 
 properties of elements of the above family.  

\begin{theorem}(\cite[Th.~3]{Zd})  \label{ch1}
Let $X$ be a Lindel\"of topological space. Then $X$ has the property $\uop$  if and only if
for every  $u\in\Lambda(X)$  so does the semifilter $\uparrow\mu_u(X)$.
\end{theorem}

And finally, we define some  properties  of semifilters closely
related to  $\utas$ and $\utass$. We say that a family $\B\subset\F$ is  a \emph{base} of
a semifilter $\F$, if $\F=\uparrow\B$.
The \emph{character} $\chi(\F)$ of a semifilter $\F$ 
equals, by definition, to the smallest size of a base of $\F$.

\begin{definition} A filter $\F$ on  $\w$ is defined to be
a \emph{simple $P$-filter}, if there exists a linearly preordered with respect
to $\subset^\ast$ base of $\F$.
\end{definition}

The subsequent observation explains the importance of simple $P$-filters
in studying of properties $\utas$ and $\utass$.

\begin{observation} \label{obs1}
A family $u=\{U_n:n\in\w\}$ of subsets of $X$ is a $\tau^\ast$- (resp. $\tau^\star$-) cover of $X$
if and only if $\mu_u(X)$ can be enlarged to (resp. generates)
 a simple $P$-filter.
\end{observation}

We shall also use the subsequent characterization of simple $P$-filters.
\begin{theorem}(\cite[3.2.3]{BZ}) \label{spf-ch}
A filter $\F$  is a simple $P$-filter if and only if
 $\F$ has a base $\B=(B_\alpha)_{\alpha<\chi(\F)}$
such that $B_\alpha\subset^\ast B_\beta$ for all
$\beta\leq\alpha<\chi(\F)$.
\end{theorem}

Next, we shall search for conditions when there are nonmeager simple
$P$-filter, or conditions which imply that all of them are meager.

\begin{proposition} \label{search-m}
If $\dep\leq\mb$, then each simple $P$-filter is meager.
\end{proposition}
\begin{proof}
Easily follows from Theorem~\ref{spf-ch}, the definition of the cardinal 
$\dep$, and the fact that each semifilter wich character $<\mb$ is meager,
see \cite[8.3.1]{Zd} or \cite{Solomon}. 
\end{proof}

\begin{proposition} \label{search-nm}
There exists a nonmeager simple $P$-filter provided $\mb<\md$ or $\mt=\mb$.
\end{proposition}
\begin{proof}
Follows immediately from  \cite[8.3.2, 11.2.3]{BZ}.
\end{proof}

The following simple characterization of the property $\uga$
is of crucial importance for the proof of Theorem~\ref{main}(3).
Let $u$ be a cover of a set $X$. A subset $B$ of $X$ is \emph{$u$-bounded},
if $B\subset\cup v$ for some finite $v\subset u$.

\begin{proposition} \label{s-c}
A topological space $X$ has  the property $\uga$ if and only if for every sequence
$(u_n)_{n\in\w}$ of open covers of $X$ there exists a sequence
$(v_n)_{n\in\w}$ such that each $v_n$ is a finite subset of $u_n$
and a semifilter $\uparrow\mu_{\{\cup v_n:n\in\w\}}(X)$ is meager.
\end{proposition}
\begin{proof}
Only the ``if'' part needs a proof.
Let  $(u_n)_{n\in\w}$
be a sequence of open covers of $X$.
Without loss of generality, $u_{n+1}$ is a refinement of $u_n $ for all $n\in\w$.
Let $w=\{B_n:n\in\w\}$ be such that
 each $B_n$ is $u_n$-bounded and $\uparrow\mu_{w}(X)$ is meager. Then there is an increasing
number sequence $(k_n)_{n\in\w}$ such that each element of $\uparrow\mu_{w}(X)$
meets all but finitely many half-intervals $[k_n,k_{n+1})$.
Since $u_{n+1}$ is a refinement of $u_n$ for all $n\in\w$,
the union $C_n=\bigcup_{k\in[k_n,k_{n+1})}B_k$ is $u_n$-bounded.
We claim that  $\{C_n:n\in\w\}$ is a $\gamma$-cover of $X$. 
Indeed, given any $x\in X$ find $n_0\in\w$ such that $\mu_{w}(x)\cap [k_n,k_{n+1})\neq\emptyset$
for all $n\geq n_0$. The above means that for every $n\geq n_0$ we can find $k_x(n)\in [k_n,k_{n+1})$
with the property $x\in B_{k_x(n)}$, and hence $x\in B_{k_x(n)}\subset\bigcup_{k\in [k_n,k_{n+1})}B_k=C_n$
for all $n\geq n_0$. 
\end{proof}

In  the proofs of Theorem~\ref{main}
we shall use some properties of the \emph{eventual dominance relation}
$\leq^\ast$ on $\IN^\w$ defined as follows:
$x\leq^\ast y$ whenever the set $\{n\in\w:x_n>y_n\}$ is finite.
A subset $A$ of $\IN^\w$ is said to be
\begin{itemize}
\item \emph{bounded}, if there exists
$x\in\IN^\w$ such that $a\leq^\ast x$ for every  $a\in A$;
\item \emph{dominating}, if 
 for every $x\in\IN^\w$  there exists $a\in A$ such that
 $x\leq^\ast a$;
\item  a \emph{scale}, if there exists an ordinal
$\alpha$ and a bijection $\varphi: \alpha\to A$ such that $\varphi(\beta)\leq^\ast\varphi(\eta)$ 
for all $\beta<\eta$.
In case $\alpha=\mb$ the set $A$ is said to be a \emph{$\mb$-scale}.
\end{itemize}
%
%
%As it was shown  by P.~L.~Dordal in \cite{dordal},
% it is consistent that $\mathfrak h=\mb=\aleph_2$ and there are no towers of length
%$\aleph_2$. Therefore it is
%consistent that there is no tower of lenght $\mb$,
%which yields the meagerness of all $\tau^\ast$-semifilters by Lemma~\ref{l2}.

%We stress that in every such a model $\mb=\md$ by Lemma~\ref{l1}.
%It is also worth to mention here that A.~Dow, extending the above result of Dordal, 
%showed in \cite{Dow} that this  holds in the Mathias model, see \cite{BlassHBK}
%and references there in.
%

\noindent \textbf{Proof of Teorem~\ref{main}.}
 Let $X$ be a topological space and $(u_n)_{n\in\w}$
be a sequence of open covers of $X$ such that $u_{n+1}$
is an refinement of $u_n$ for all $n\in\w$. 
\smallskip

 1. As it was mentioned in Intoduction, $(\muu<\mg)$ implies $(\mb<\md)$,
and therefore there exists a nonmeager simple $P$-filter
$\F$ by Proposition~\ref{search-nm}. 
By the definition of the property $\uop$
there exists and a large cover $w_1=\{B_n:n\in\w\}$ of $X$
such that each $B_n$ is $u_n$-bounded, see \cite{Sc}.
 Applying Theorem~\ref{ch1} we conclude that the semifilter
$\U=\uparrow\mu_{w_1}(X)$ has the property $\uop$, and consequently
it is not comeager by \cite[Prop.~2]{Zd}. 
Two cases are possible.

a) $\U$ Is bi-Baire. Then Theorem~\ref{lll} supplies us with an
increasing sequence $(k_n)_{n\in\w}$ such that 
$\G:=\phi(\U)=\phi(\F)$, where
$\phi:\w\to\w$ is such that $\phi^{-1}(n)=[k_n,k_{n+1})$ for all $n\in\w$,
and $\phi(\A)=\{\phi(A):A\in\A\}$ for any family $\A$ of subsets of $\w$. 
Note that $\G$ is  a simple $P$-filter being an image of $\F$
under $\phi$. 

Let $C_n=\bigcup_{k\in [k_n,k_{n+1})}B_k$.
By our choice of $(u_n)_{n\in\w}$, each $C_n$ is $u_n$-bounded. 
We claim that $w_2=\{C_n:n\in\w\}$ is a $\tau^\star$-cover of $X$. 
Indeed, since $\G=\phi(\U)$, $\U$ is generated by
$\mu_{w_1}(X)$, and
$\mu_{w_2}(x)=\phi(\mu_{w_1}(x))$ for all $x\in X$,
we conclude that $\G$ is generated by $\mu_{w_2}(X)$.  
Now it sufficies to apply Observation~\ref{obs1}. 
 
b) $\uparrow\mu_{w_1}(X)$ is meager.
Then in the same way as in the proof of Proposition~\ref{s-c}
we can construct a $\gamma$-cover $\{C_n:n\in\w\}$ of $X$
such that each $C_n$ is $u_n$-bounded.
\smallskip

2. In this case  it sufficies to find an $\w$-cover $w_1=\{B_n:n\in\w\}$ of $X$
such that each $B_n$ is $u_n$-bounded and apply to the filter
$\uparrow\mu_{w_1}(X)$ the same arguments as in the proof of the first item.  
\smallskip

3. Let us assume that each simple $P$-filter is meager and $X$
has the property $\utas$.
Then there exists a $\tau^\ast$-cover $w=\{B_n:n\in\w\}$ of $X$ such that each $B_n$
is $u_n$-bounded. 
By Observation~\ref{obs1} this implies  that the semifilter $\U=\uparrow\mu_w(X)$
can be enlarged to  a simple $P$-filter $\F$, which is meager by our assumption,
 and hence so is $\U$.
Applying  Proposition~\ref{s-c} we conclude that $X$ has the property $\uga$.

Next, suppose that there exists a nonmeager simple $P$-filter $\F$.
The rest of the proof falls naturally into two parts.
\smallskip

a) $(\mb=\md).$  In this case the assertion follows from \cite[8.10]{TZ},
which supplies us with a subspace
$Y$ of the Baire space with the following properties:
\begin{itemize}
\item[$(i)$] $Y$ does not have the property $\uta$;
\item[$(ii)$] for any sequence $(w_n)_{n\in\w}$ of open covers 
of $Y$ there exists a family $w=\{B_n:n\in\w\}$ such that
each $B_n$ is $w_n$-bounded and
$\uparrow\mu_w(X)\subset\F$.
\end{itemize}

b) $(\mb<\md).$  In this case the assertion follows from
the subsequent two statements.
\begin{itemize}
\item[$(i)$] There exists a subspace of the Baire space of size $\mb$
which does not have the property $\uga$.
\item[$(ii)$] $(\mb<\md)$ Implies that every subspace $Y$ of the Baire
space satisfies $\utas$ provided $|Y|\leq \mb$.
\end{itemize}
The first of them may be found in \cite{Sc}. To prove the second one,
find a (probably not bijective) enumeration $\{y_\alpha:\alpha<\mb\}$
of $Y$. Recall from \cite{Ts-tow} that a subset $Z\subset\w^\w$
has a \emph{weak excluded middle property} if there exists 
$x\in \w^\w$ such that the family $\{[z\leq x]:z\in Z\}$ can be enlarged
to a simple $P$-filter, where for a relation $R$ on $\w$
$[z\: R\: x]=\{n\in\w:z(n)\: R\: x(n)\}$.

Let $f:Y\to\w^\w$ be continuous. By transfinite induction over $\mb$
construct a  $\mb$-scale $B=\{b_\alpha:\alpha<\mb\}$ such that
$f(y_\alpha), b_\beta\leq^\ast b_\alpha$ 
for all $\beta\leq\alpha<\mb$.
Since $\mb<\md$, $B$ is not dominating, which means that there exists
$c\in\w^\w$ such that $c\leq^\ast b_\alpha$ for no $\alpha<\mb$,
and hence $[b_{\alpha}< c]$ is infinite
for all $\alpha$. Observe that for arbitrary $\beta\leq\alpha<\mb$
the equation  $b_\beta\leq^\ast b_\alpha$ implies $[b_\alpha<c]\subset^\ast[b_\beta<c]$,
and therefore $\mathcal T=([b_\alpha<c])_{\alpha<\mb}$ is a tower.
Moreover, $[b_\alpha<c]\subset^\ast[f(y_\alpha)\leq c]$,
 consequently the family $\{[f(y_\alpha)\leq c]:\alpha<\mb\}=\{[f(y)\leq c]:y\in Y\}$
is a subset of the simple $P$-filter generated by $\mathcal T$, and
hence $f(Y)$ has a weak excluded middle property. 
Applying \cite[Th.~7.8]{Ts-tow} asserting
that a subset $Z$ of the Baire space satisfies $\utas$ provided
for every continuous $\phi:Z\to\w^\w$  the image $\phi(Z)$ has the weak
excluded middle property, we conclude that $Y$ has the property $\utas$.
\hfill $\Box$
\medskip

\noindent \textbf{Proof of Corollary~\ref{details}.}

1. Follows immediately from Proposition~\ref{search-m} and
Theorem~\ref{main}(3).
\smallskip

2. Under $(\mb<\md)$ the assertion follows from
Proposition~\ref{search-nm} and Theorem~\ref{main}(3).

Under $(\mt=\md)$ it sufficies to use $(\mt=\mb)$--part of Proposition~\ref{search-nm}
to find a nonmeager simple $P$-filter and then apply the same
arguments as in the proof of the $(\mb=\md)$--part of Theorem~\ref{main}(3). 
\hfill $\Box$
\medskip

\noindent \textbf{Acknowledgements.} The author wishes to express his thanks to prof.
Taras Banakh for supervising the writing of this paper, and
  prof. Boaz Tsaban
 for  fruitful discussions and suggestions.

\noindent Department of Mechanics and Mathematics,\\
 Ivan Franko Lviv National University, \\
Universytetska 1, Lviv, 79000, Ukraine.
\medskip

\textit{E-mail address:}   \texttt{lzdomsky@rambler.ru}

\footnotesize 

\begin{thebibliography}{ChGP??}

\bibitem[1]{BZ} Banakh~T., Zdomsky~L., {\it Coherence of semifilters},\\ 
\texttt{http://www.franko.lviv.ua/faculty/mechmat/Departments/Topology/booksite.html} 


\bibitem[2]{BZ-mult} Banakh~T., Zdomsky~L., {\it
Selection principles and infinite games on multicovered spaces 
and their applications,} in preparation.

\bibitem[3]{BST} Bartoszy\'nski~T., Shelah~S., Tsaban~B., {\it Additivity properties of topological diagonalizations},
 J. Symbolic Logic \textbf{68} (2003), 1254--1260. \\
(Full version: \texttt{http://arxiv.org/abs/math.LO/0112262})

\bibitem[4]{BlassHBK}  Blass~A., 
  \emph{Combinatorial cardinal characteristics of the continuum},
  in Handbook of Set Theory (M.~Foreman et. al., Eds.),
    to appear.

%\bibitem[5]{BlS} A.~Blass, S.~Shelah, {\it There may be simple $P_{\aleph_1}$- and $P_{\aleph_2}$-points
%and the Rudin-Keisler ordering may be downward directed,} Ann. Pure Appl. Logic
%\textbf{33} (1987), 213--243.

\bibitem[5]{CP} Chaber~J., Pol~R., {\it A remark on Fremlin-Miller theorem
concerning the Menger property and Michael concentrated sets}, preprint.

\bibitem[6]{dordal}  Dordal P., %no tower of size h=aleph_2 in Mathias
 {\it A model in which the base-matrix tree cannot have cofinal branches},
  J. Symbolic Logic \textbf{52} (1987), 651--664.

\bibitem[7]{Dow} %no tower of size h=aleph_2 in Mathias
  Dow A., \emph{Set theory in topology},
  in Recent Progress in General Topology (M.~Husek et. al., Eds.),
    Elsevier Sci. Publ., Amsterdam, 1992, pp. 168--197.

\bibitem[8]{FM} Miller~A., Fremlin~D., {\it On some properties of Hurewicz,
    Menger, and Rothberger\/}, Fund. Math. \textbf{129} (1988), 17--33.

\bibitem[9]{GN} Gerlits~J., Nagy~Zs., {\it
Some properties of $C(X)$, I},  
Topology Appl. $\mathbf{14}$ (2) (1982), 151--163.  

\bibitem[10]{Hu} Hurewicz~W., {\it \"{U}ber die Verallgemeinerung des Borelschen
Theorems}, Mathematische Zeitschrift \textbf{24} (1925) 401-421.

\bibitem[11]{JMSS} Just~W., Miller~A., Scheepers~M., Szeptycki~S.,
        {\it The combinatorics of open covers II\/}, Topology  Appl.
       \textbf{73} (1996), 241--266.

\bibitem[12]{La} Laflamme~C., {\it Equivalence of families of
functions on natural numbers}, Trans. Amer. Math. Soc. \textbf{330} (1992),
307--319.

\bibitem[13]{M-S} Marczewski~E. (Szpilrajn), {\it The characteristic function
of a sequence of sets and some of its applications,} Fund. Math.
\textbf{31} (1938), 207--233.


\bibitem[14]{Me} Menger~K., {\it Einige \"{U}berdeckungss\"{a}tze der
Punktmengenlehre, Sitzungsberichte}. Abt. 2a, Mathematic,
Astronomie,
 Physic, Meteorologie und Mechanic (Wiener Akademie) \textbf{133} (1924) 421-444.


\bibitem[15]{Sc} Scheepers~M., {\it Combinatorics of open covers I: Ramsey Theory},
Topology Appl. \textbf{69}~(1996), 31--62.

\bibitem[16]{ST} Shelah~S., Tsaban~B., {\it Critical cardinalities and additivity
properties of combinatorial notions of smallness}, J. Appl. Anal.
\textbf{9} (2003), 149--162.\\
\texttt{http://arxiv.org/abs/math.LO/0304019}

\bibitem[17]{Solomon} Solomon~R., {\it Families of sets and functions,}
Czechoslovak Math. J. \textbf{27}~(1977), 556-559. 

\bibitem[18]{Ta} Talagrand~M., {\it Filtres: Mesurabilit\'{e}, rapidit\'{e},
 propri\'{e}t\'{e} de Baire forte}, Studia Math. \textbf{74} (1982), 283--291.

\bibitem[19]{Ts-tow} Tsaban~B., {\it Selection principles and the minimal tower problem,}
Note Math., to appear.  \texttt{http://arxiv.org/abs/math.LO/0105045}

\bibitem[20]{SPM3} Tsaban~B., (eds.), \\
SPM Bulletin \textbf{3} (2003) \texttt{http://arxiv.org/abs/math.GN/0303057}

\bibitem[21]{TZ} Tsaban~B., Zdomsky~L., {\it Scales, Fields, and a problem of Hurewicz},
submitted to J. Amer. Math. Soc.. \texttt{http://arxiv.org/abs/math.GN/0507043}. 

\bibitem[22]{Va} Vaughan J.,
  {\it Small uncountable cardinals and topology},
  in Open problems in topology (J. van Mill, G.M. Reed,
Eds.),
    Elsevier Sci. Publ., Amsterdam, 1990, pp. 195-218.


\bibitem[23]{Zd} Zdomsky~L., {\it A semifilter approach to selection principles,} 
to appear in Comment. Math. Univ. Carolinae.  
\texttt{http://arxiv.org/abs/math.GN/0412498}

\end{thebibliography}
\end{document}